\documentclass[12pt,reqno,a4wide]{amsart}

\usepackage{tikz} 
\usepackage{calrsfs}
\usepackage[charter]{mathdesign}

\allowdisplaybreaks

\title{Representing derivatives of Chebyshev polynomials by Chebyshev polynomials}

\author[H.~Prodinger]{Helmut Prodinger}
\address{H.~Prodinger\\Department of Mathematical Sciences, Mathematics Division\\ Stellenbosch University, Private Bag X1, 7602 Matieland, South Africa}
\email{hproding@sun.ac.za}

\thanks{}
\date{\today}

\begin{document}
	\begin{abstract}
		A recursion formula for derivatives of Chebyshev polynomials is replaced by an  explicit formula.
	\end{abstract}
	
	\maketitle
	
	\section{Introduction}
	
	Consider the Chebyshev polynomials of the second kind
\begin{equation*}
U_n(x)=\sum_{0\le k\le n/2}(-1)^k\binom{n-k}{k}(2x)^{n-2k};
\end{equation*}
the main interest of the paper~\cite{Siyi15} is to represent the derivatives of $U_n(x)$ in terms of the Chebyshev polynomials themselves. To this aim an ``exact computational method'' (a recursion formula) was presented. In the present note, we replace this ``computational method'' by an exact and explicit formula.

Our answer is
\begin{align*}
U_n^{(s)}(x)
&=2^{s}\sum_{0\le j \le (n-s)/2}(n-j)^{\underline{s-1}}\binom{s+j-1}{s-1}(n-2j-s+1)
U_{n-s-2j}(x).
\end{align*}

Although it is not needed, we briefly mention without  proof an analogous formula for the Chebyshev polynomials of the first kind: Let
\begin{equation*}
T_n(x)=\sum_{0\le k\le n/2}(-1)^k\frac{n}{n-k}\binom{n-k}{k}2^{n-1-2k}x^{n-2k},
\end{equation*}
then
\begin{align*}
T_n^{(s)}(x)&=2^s\sum_{0\le j\le (n-s)/2}n(n-1-j)^{\underline{s-1}}\binom{s+j-1}{s-1}T_{n-s-2j}(x)\\
&-[\![n-s\text{ even}]\!]\;2^{s-1}n\,((n+s)/2-1)^{\underline{s-1}}\binom{(n+s)/2-1}{s-1}.
\end{align*}

We use here the notion of falling factorials $x^{\underline n}:=x(x-1)\dots (x-n+1)$ and Iverson's symbol
$[\![P]\!]$ which is 1 if $P$ is true and 0 otherwise, compare \cite{GrKnPa94}.

\section{The proof}

Our starting point is the inversion formula (see \cite{Riordan68})
\begin{equation*}
x^j=2^{-j}\sum_{0\le h \le j/2}\bigg[\binom{j}{h}-\binom{j}{h-1}\bigg]U_{j-2h}(x),
\end{equation*}
which we will use in
\begin{equation*}
U_n^{(s)}(x)=\sum_{0\le k\le n/2}(-1)^k\binom{n-k}{k}(n-2k)^{\underline{s}}2^{n-2k}x^{n-s-2k}
\end{equation*}
and simplify:

\begin{align}\label{arsch}
U_n^{(s)}(x)
&=\sum_{0\le k\le n/2}(-1)^k\binom{n-k}{k}(n-2k)^{\underline{s}}
2^{s}\sum_{k\le h+k \le (n-s)/2}\bigg[\binom{n-s-2k}{h}-\binom{n-s-2k}{h-1}\bigg]U_{n-s-2k-2h}(x)\notag\\
&=\sum_{0\le k\le n/2}(-1)^k\binom{n-k}{k}(n-2k)^{\underline{s}}
2^{s}\sum_{k\le j \le (n-s)/2}\bigg[\binom{n-s-2k}{j-k}-\binom{n-s-2k}{j-k-1}\bigg]U_{n-s-2j}(x)\notag\\
&=2^{s}\sum_{0\le k\le j \le (n-s)/2}(-1)^k\binom{n-k}{k}(n-2k)^{\underline{s}}
\bigg[\binom{n-s-2k}{j-k}-\binom{n-s-2k}{j-k-1}\bigg]U_{n-s-2j}(x).
\end{align}
We compute the sum over $k$ separately:
\begin{align*}
\sum_{0\le k\le j }&(-1)^k\binom{n-k}{k}(n-2k)^{\underline{s}}
\bigg[\binom{n-s-2k}{j-k}-\binom{n-s-2k}{j-k-1}\bigg]\\
&=\sum_{0\le k\le j }(-1)^k
\bigg[\frac{(n-k)!}{k!(j-k)!(n-s-k-j)!}-\frac{(n-k)!}{k!(j-1-k)!(n-s-k-j+1)!}\bigg]\\
&=\sum_{0\le k\le j }(-1)^k
\bigg[\frac{(n-j)!}{(n-s-j)!}\binom{n-k}{j-k}\binom{n-s-j}{k}
-\frac{(n-j+1)!}{(n-s-j+1)!}\binom{n-k}{j-1-k}\binom{n-s-j+1}{k}\bigg]\\
&=\sum_{0\le k\le j }(-1)^j
\bigg[(n-j)^{\underline{s}}\binom{-n+j-1}{j-k}\binom{n-s-j}{k}
+(n-j+1)^{\underline{s}}\binom{-n+j-2}{j-1-k}\binom{n-s-j+1}{k}\bigg]\\
&=(-1)^j\bigg[(n-j)^{\underline{s}}\binom{-1-s}{j}
+(n-j+1)^{\underline{s}}\binom{-1-s}{j-1}\bigg]\\
&=(n-j)^{\underline{s}}\binom{s+j}{j}
-(n-j+1)^{\underline{s}}\binom{s+j-1}{j-1}\\
&=(n-j)^{\underline{s-1}}\binom{s+j-1}{s-1}\bigg[(n-j-s+1)\frac{s+j}{s}-(n-j+1)\frac{j}{s}\bigg]\\
&=(n-j)^{\underline{s-1}}\binom{s+j-1}{s-1}(n-2j-s+1).
\end{align*}
In this computation only the Vandermonde convolution formula \cite{GrKnPa94} was used.

Plugging this formula into (\ref{arsch}) yields the announced formula from the introduction.

\end{document}